\documentclass[11pt, a4paper]{article}

\usepackage[top=2.6cm, bottom=2.6cm, left=2.5cm, right=2.5cm]{geometry}

\usepackage{latexsym,amssymb}
\usepackage{amsmath,amsthm}
\usepackage{color}
\usepackage{mathtools}

\usepackage{multirow,eepic}

\newtheorem{Theorem}{\bf Theorem}[section]

\newtheorem{Lemma}{\bf Lemma}[section]
\newtheorem{Proposition}{\bf Proposition}[section]
\newtheorem{Corollary}{\bf Corollary}[section]
\newtheorem{Remark}{\bf Remark}[section]

\allowdisplaybreaks[4]

\numberwithin{equation}{section}
\begin{document}
\title{
  Quasi self-similarity and its application to 
  the global in time solvability of a superlinear heat equation
}
\author{
        Yohei Fujishima\footnote{e-mail address: fujishima@shizuoka.ac.jp} \\ \\
        {\small Department of Mathematical and Systems Engineering} \\
        {\small Faculty of Engineering, Shizuoka University} \\
        {\small 3-5-1 Johoku, Hamamatsu 432-8561, Japan} \\ \\
        Norisuke Ioku\footnote{e-mail address: ioku@tohoku.ac.jp} \\ \\
        {\small Mathematical Institute, Tohoku University} \\
        {\small Aramaki 6-3, Sendai 980-8578, Japan}
        }
\date{}
\pagestyle{myheadings}
\markboth{Y. Fujishima and N. Ioku}{Quasi self-similarity for nonlinear heat equations}

\maketitle

\begin{abstract}
  This paper concerns the global in time existence of solutions for 
  a semilinear heat equation 
  \begin{equation}
    \tag{P}
    \label{eq:P}
    \begin{cases}
      \partial_t u = \Delta u + f(u),
      &x\in \mathbb{R}^N, \,\,\, t>0,
      \\[3pt]
      u(x,0) = u_0(x) \ge 0,
      &x\in \mathbb{R}^N,
    \end{cases}
  \end{equation}
  where $N\ge 1$, $u_0$ is a nonnegative initial function
  and $f\in C^1([0,\infty)) \cap C^2((0,\infty))$ denotes superlinear nonlinearity of the problem. 
  We consider the global in time existence and nonexistence of solutions for problem~\eqref{eq:P}. 
  The main purpose of this paper is to determine 
  the critical decay rate of initial functions for the global existence of solutions. 
  In particular, we show that it is characterized by 
  quasi self-similar solutions which are solutions $W$ of 
  \begin{equation}
    \notag
    \Delta W + \frac{y}{2}\cdot \nabla W + f(W)F(W) + f(W)
    + \frac{|\nabla W|^2}{f(W)F(W)}
    \Bigl[
      q - f'(W)F(W)
    \Bigr]
    = 0, 
    \quad 
    y \in \mathbb{R}^N, 
  \end{equation}
  where 
  $F(s):=\displaystyle\int_s^{\infty}\dfrac{1}{f(\eta)}d\eta$
  and
  $q\ge 1$. 
\end{abstract}
\noindent
{\bf Keywords}: nonlinear heat equation,
quasi-scaling, forward self-similar solutions, 
global existence of solutions, 
blow-up of solutions 
\vspace{5pt}
\newline
{\bf 2010 MSC}: Primary; 35K91, Secondly; 35B44, 35C06

\ \\
Correspondence and proofs to: Yohei Fujishima $<$fujishima@shizuoka.ac.jp$>$

\section{Introduction}
Consider
\begin{equation}
  \label{eq:1.1}
  \begin{cases}
    \partial_t u = \Delta u + f(u),
    &x\in \mathbb{R}^N, \,\,\, t>0,
    \\[3pt]
    u(x,0) = u_0(x) \ge 0,
    &x\in \mathbb{R}^N,
  \end{cases}
\end{equation}
where $N\ge 1$, $u_0$ is a nonnegative initial function
and $f\in C^1([0,\infty)) \cap C^2((0,\infty))$ is a function
denoting the nonlinearity of problem~\eqref{eq:1.1}.
Throughout this paper, we assume that $f$ satisfies 
\begin{equation}
  \label{eq:1.2}
  f(0) = 0 
  \quad\mbox{and}\quad 
  f(u)>0,
  \quad
  f'(u)>0,
  \quad
  f''(u)>0,
  \quad \mbox{for all}\,\,\, u>0.
\end{equation}
Furthermore, we assume that $f$ satisfies
\begin{equation}
  \label{eq:1.3}
  F(u) := \int_u^\infty \frac{1}{f(s)} \, ds < \infty
  \quad\mbox{for all}\,\,\, u>0.
\end{equation}
By \eqref{eq:1.2} and \eqref{eq:1.3} we have $F'(u) = -1/f(u)<0$.
Moreover, 
since $F(u) \to \infty$ as $u\to +0$ by \eqref{eq:1.2} (see \cite[Lemma~3.1]{FI3})
and $F(u) \to 0$ as $u\to \infty$,
the range of $F$ is $(0,\infty)$.
So the function $F$ is monotonically 
decreasing
and its inverse function $(0,\infty) \ni v \mapsto F^{-1}(v) \in (0,\infty)$,
which is also decreasing function, exists.
This paper is devoted to the study of the relationship between
the global in time solvability for problem~\eqref{eq:1.1} and
the decay rates of initial functions at space infinity.
In particular, we characterize the critical decay rate of initial data 
for the global in time existence of solutions for problem~\eqref{eq:1.1} 
by using the function $F^{-1}$. 

This paper concerns the global in time existence of solutions for problem~\eqref{eq:1.1}, 
which has been studied intensively for various nonlinear terms such as 
power nonlinearity (\cite{BP}, \cite{Fujita}, 
\cite{HW}, \cite{H}, \cite{KST}, \cite{N}, \cite{N2}, \cite{Sugitani} and \cite{W2}), 
nonlinearity with a gradient term (\cite{STW} and \cite{Ta}), 
exponential nonlinearity (\cite{F2}, \cite{I}, \cite{IRT} and \cite{V}) 
and general nonlinearity (\cite{FI3} and \cite{LS2}). 
See also the book \cite{QS}, which includes further information and a list of numerous references on the topic. 
Among others, the global in time existence of solutions for a semilinear heat equation is 
mainly studied for the case of $f(u) = u^p$ with $p>1$,
that is, 
\begin{equation}
  \label{eq:1.4}
  \begin{cases}
    \partial_t u = \Delta u + u^p,
    &x\in \mathbb{R}^N, \,\,\, t>0,
    \\[3pt]
    u(x,0) = u_0(x) \ge 0,
    &x\in \mathbb{R}^N,
  \end{cases}
\end{equation}
where $p>1$. 
Fujita in \cite{Fujita} showed that the critical exponent $p_F := 1+2/N$ classifies the 
existence of global solutions for a power type semilinear heat equation. 
See also \cite{H}, \cite{KST} and \cite{Sugitani}. 
In particular, if $p>p_F$, then there exists a global solution of \eqref{eq:1.4} 
for rapidly decaying initial data. 
Weissler in \cite{W2} gave a characterization of a spatial decay rate for initial data in terms of the integrability. 
It has been shown that solutions exist globally in time for small initial data in the Lebesgue space 
$L^{\frac{N}{2}(p-1)}(\mathbb{R}^N)$.
A functional space $L^{\frac{N}{2}(p-1)}(\mathbb{R}^N)$, 
which is invariant under a self-similar scaling (see below), 
is also known as the critical space for the local in time well-posedness (see \cite{W1}). 
Lee--Ni in~\cite{LN} specified the optimal decay rate of initial data for the global in time existence of solutions, 
which is given by $|x|^{-\frac{2}{p-1}}$. 
Namely, 
for an initial function $u_0$ which satisfies $u_0(x) \simeq \gamma |x|^{-\frac{2}{p-1}}$ 
as $|x| \to \infty$ for a constant $\gamma>0$, 
there exists a global in solution of \eqref{eq:1.4} with small enough $\gamma$,  
while it blows up in finite time if $\gamma$ is large enough. 
Moreover, the optimal constant $\gamma$ classifying the existence of global solutions was 
characterized by Naito~\cite{N2}. 

We remark that the optimal decay rate of initial data 
is based on
the scaling property for a power type semilinear heat equation. 
Let $u$ satisfy
\begin{equation*}
  \partial_t u = \Delta u + u^p \qquad \mbox{in} \quad \mathbb{R}^N \times (0,\infty).  
\end{equation*}
A scaled function $u_\lambda(x,t):= \lambda^\frac{2}{p-1}u(\lambda x,\lambda^2 t)$ with $\lambda>0$, 
which is called a self-similar scaling, 
satisfies the same equation, 
and its $L^r(\mathbb{R}^N)$ norm at $t=0$ is invariant if and only if $r=\frac{N}{2}(p-1)$ when $p>p_F$.  
Then its dilation rate $\frac{2}{p-1}$ gives the optimal decay rate for initial data as is shown by Lee--Ni~\cite{LN}. 
The optimal
decay rate classifying the global existence of solutions for problem~\eqref{eq:1.1} with $f(u) = e^u$ has been studied by 
the first author of this paper in \cite{F2} 
and it is given by $-2\log |x|$, 
where the equation also possesses a self-similar property. 
Consequently, the optimal decay rate of $u_0$ for the global existence of solutions of 
problem~\eqref{eq:1.1} with $f(u) = u^p$ or $f(u) = e^u$ 
is characterized by a self-similar property of the equation, 
and it is also given by the decay rate of forward self-similar solutions.  
Forward self-similar solutions are functions satisfying 
\begin{equation}
  \label{eq:1.5}
  \Delta W + \frac{y}{2}\cdot\nabla W + \frac{1}{p-1} W + W^p = 0,
  \qquad y\in \mathbb{R}^N, 
\end{equation}
if $f(u)=u^p$ 
and 
\begin{equation}
  \label{eq:1.6}
  \Delta W + \frac{y}{2}\cdot\nabla W + 1 + e^W = 0,
  \qquad y\in \mathbb{R}^N, 
\end{equation}
if $f(u)=e^u$. 
The behavior of forward self-similar solutions as $|x| \to \infty$ is given by $|x|^{-\frac{2}{p-1}}$ if $f(u) = u^p$ 
and by $-2\log |x|$ if $f(u) = e^u$. 

Concerning problem~\eqref{eq:1.1}, 
we can not expect 
an exact self-similar property 
except for power and exponential nonlinearities. 
This is a crucial difficulty to deal with general nonlinearities. 
On the other hand, it is possible to propose a generalization of a self-similar transformation. 
The following scaling was given in \cite{F}:
\begin{equation}
  \label{eq:1.7}
  u_\lambda(x,t):= 
  F^{-1} \Bigl[
    \lambda^{-2} F(
      u(\lambda x,\lambda^2 t)
    )
  \Bigr], 
  \qquad 
  \lambda>0, 
\end{equation}
which satisfies 
\begin{equation*}
  \partial_t u_\lambda = \Delta u_\lambda + f(u_\lambda) + 
  \frac{|\nabla u_\lambda|^2}{
    f(u_\lambda)F(u_\lambda)
  } \Bigl[
    f'(u)F(u) - f'(u_\lambda)F(u_\lambda)
  \Bigr]. 
\end{equation*}
A method to construct a supersolution based on the transformation~\eqref{eq:1.7} was 
proposed by the authors of this paper in \cite{FI}, 
and the local and global in time existence of solutions for problem~\eqref{eq:1.1} 
is studied in \cite{FI} and \cite{FI3}.
In particular, the authors of this paper proved 
in \cite{FI3} that 
a generalization of the Fujita exponent for problem~\eqref{eq:1.1} 
is given by 
\begin{equation}
  \label{eq:1.8}
  \lim_{s\to +0} f'(s) F(s) = \frac{p_F}{p_F-1} = 1+\frac{N}{2}, 
\end{equation}
and there exists a global in time solution of \eqref{eq:1.1} if 
\begin{equation}
  \label{eq:1.9}
  \lim_{s\to +0} f'(s) F(s) < 1+\frac{N}{2}. 
\end{equation}
The critical case~\eqref{eq:1.8} still contains open problems and partial results are known (see \cite[Theorem~1.3]{FI3}). 
Other applications of the transformation \eqref{eq:1.7} can be found in \cite{FI2,GM,M}. 

The purpose of this paper is to characterize the critical decay rate of initial data 
for the global existence of solution of problem~\eqref{eq:1.1} under the condition~\eqref{eq:1.9}. 
As analogies with the cases that $f(u)=u^p$ and $f(u) = e^u$, 
forward self-similar solutions seem to give the optimal decay rate of an initial function 
for the global in time existence of solutions for problem~\eqref{eq:1.1}. 
As mentioned above, 
a big issue for this conjecture is the lack of self-similarity of the equation. 
We overcome this difficulty by considering a generalization of forward self-similar solutions 
based on the transformation~\eqref{eq:1.7}. 
In fact, we prove that 
the optimal decay is determined by the behavior of solutions as $|y|\to \infty$ for an elliptic problem 
\begin{equation}
  \label{eq:1.10}
  \Delta W + \frac{y}{2}\cdot \nabla W + f(W)F(W) + f(W)
  + \frac{|\nabla W|^2}{f(W)F(W)}
  \Bigl[
    q - f'(W)F(W)
  \Bigr]
  = 0 
\end{equation}
in $y\in \mathbb{R}^N$, where $q\ge 1$. 
Equation~\eqref{eq:1.10} is equal to \eqref{eq:1.5} and \eqref{eq:1.6} 
if $f(u)=u^p$ with $q=\frac{p}{p-1}$ and $f(u)=e^u$ with $q=1$, respectively. 
From this point of view, 
we call a solution of \eqref{eq:1.10} a \textit{quasi forward self-similar solution} for problem~\eqref{eq:1.1}.

\medskip 

We are ready to state main results of this paper. 
In what follows, we assume that there exists a limit
\begin{equation}
  \label{eq:1.11}
  q:= \lim_{s\to +0} f'(s) F(s).
\end{equation}
If the limit \eqref{eq:1.11} exists,
then it follows that $q \ge 1$
(see Lemma~\ref{Lemma:2.1}).
We first introduce a result on the global in time existence of solutions of problem~\eqref{eq:1.1}. 

\begin{Theorem}
  \label{Theorem:1.1}
  Let $N\ge 1$ and $f\in C^1([0,\infty)) \cap C^2((0,\infty))$ satisfy \eqref{eq:1.2} and \eqref{eq:1.3}.
  Assume that the limit 
  $q$ defined by \eqref{eq:1.11} exists
  and it satisfies
  \begin{equation}
    \label{eq:1.12}
    q < 1+\frac{N}{2}.
  \end{equation}
  For $\gamma>0$, let 
  \begin{equation}
    \notag 
    u_0(x) = 
    F^{-1} \Bigl[
      \gamma^{-1} (|x|^2+1)
    \Bigr]
    \quad\mbox{in}\quad 
    \mathbb{R}^N. 
  \end{equation}
  Then 
  there exists a constant $\gamma_*>0$ such that 
  the solution of problem~\eqref{eq:1.1} exists globally in time for $\gamma\in (0,\gamma_*)$. 
\end{Theorem}

We can also obtain behavior
of the solution constructed in Theorem~\ref{Theorem:1.1}.
In particular, the solution $u$ is bounded above and below by
quasi forward self-similar solutions characterized by \eqref{eq:1.10}.

\begin{Theorem}
  \label{Theorem:1.2}
  Assume the same conditions of Theorem~{\rm \ref{Theorem:1.1}}.
  Furthermore, assume that there exists a constant $s_0>0$ such that
  \begin{equation}
    \notag
    f'(s)F(s)\ge 1
  \end{equation}
  for $0<s<s_0$ if $q=1$.
  Then there exist
  $q^*,q_*$ with $1\le q_*\le q \le q^*<1+N/2$
  such that
  the global solution $u$ constructed in Theorem~{\rm \ref{Theorem:1.1}} satisfies
  \begin{equation}
    \label{eq:1.13}
    \frac{|y|^2}{F(W_*(y))}
    \le
    \frac{|x|^2}{F(u(x,t))}
    \le
    \frac{|y|^2}{F(W^*(y))},
    \quad
    \,\,\, y=\frac{x}{\sqrt{t+1}},
  \end{equation}
  for all $(x,t) \in \mathbb{R}^N \times [0,\infty)$,
  where $W_*$ and $W^*$ is a solution of 
  the equation \eqref{eq:1.10} with $q$ replaced by $q_*$ and $q^*$,
  respectively.
\end{Theorem}


Concerning the nonexistence of global in time solutions,
we have the following result. 

\begin{Theorem}
  \label{Theorem:1.3}
  Let $N\ge 1$ and $f\in C^1([0,\infty)) \cap C^2((0,\infty))$ satisfy \eqref{eq:1.2} and \eqref{eq:1.3}.
  Furthermore, assume that there exists a constant $s_0>0$ such that
  \begin{equation}
    \notag
    f'(s)F(s)\ge 1
  \end{equation}
  for $0<s<s_0$ if $q=1$.
  Then there exists a constant  
  $\Gamma_*$ such that,
  for all initial functions $u_0\in L^\infty(\mathbb{R}^N)$ satisfying
  \begin{equation}
    \notag    
    \liminf_{|x|\to \infty} \frac{|x|^{2}}{F(u_0(x))} \ge \Gamma_*,
  \end{equation}
  there can not exist global in time solutions for problem~\eqref{eq:1.1}.
\end{Theorem}

\begin{Remark}
  For positive constants $c_0$ and $c_1$,
  let $u_0\in L^\infty(\mathbb{R}^N)$ be a nonnegative function satisfying 
  \[
    \frac{|x|^{2}+1}{F(u_0(x)) }
    \le 
    c_0
    \iff
    u_0(x)
    \le
    F^{-1} \Bigl[
      c_0^{-1} (|x|^2+1)
    \Bigr],
  \]
  or
  \[ 
  c_1
    \le 
    \frac{|x|^{2}+1}{F(u_0(x)) }
    \iff
  F^{-1} \Bigl[
        c_1^{-1} (|x|^2+1)
      \Bigr]
    \le
    u_0(x),
  \]
  in $\mathbb{R}^N$.
  Then, combining Theorems~{\rm \ref{Theorem:1.1}} and {\rm \ref{Theorem:1.3}},
  we see that there exists a global in time solution with a sufficiently small $c_0$, 
  while the solution $u$ blows up in finite time if $c_1$ is large enough.
  Therefore, 
  the function
  $F^{-1}\left[c^{-1}(|x|^2+1)\right]$
  gives
  the critical decay rate of initial functions at space infinity for the global in time solvability for problem~\eqref{eq:1.1}.
  Since 
  \[
    F^{-1} \bigl[c^{-1} (|x|^2+1)\bigr]
    =
    \left(\frac{1}{p-1}\right)^\frac{1}{p-1}
    c^{\frac{1}{p-1}}(1+|x|^2)^{-\frac{1}{p-1}}
  \]
  if $f(u)=u^p$,
  our results can be regarded as a generalization of Lee--Ni's classification 
  for the existence of solutions for power type nonlinear heat equation. 
  The threshold constant, which is obtained by Naito~{\rm \cite{N2}}
  for $f(u)=u^p$, still remains open for a general nonlinear term $f(u)$.
\end{Remark}

As application of our results, 
we consider problem~\eqref{eq:1.1} with the following nonlinear terms: 
\begin{equation}
  \label{eq:1.14}
  f(u) = \exp\left( -\frac{1}{u} \right),
  \qquad u\in (0,1), 
\end{equation}
and 
\begin{equation}
  \label{eq:1.15}
  f(u) = u^p \left[ \log \left( e + \frac{1}{u} \right) \right]^{-r},
  \qquad u\in (0,1), 
\end{equation}
where $p>1+2/N$ and $r>0$. 
Then we see that the critical decay rate for the global in time existence of solutions for problem~\eqref{eq:1.1} 
is determined by the behavior of $f=f(u)$ as $u\to +0$. 
In particular, 
\begin{equation}
  \notag 
  u_0(x) \simeq 
  \frac{1}{\log \Big[ |x|^2 \{ \log |x| \}^2 \Big]}
  \quad\mbox{and}\quad  
  u_0(x) \simeq 
  |x|^{-\frac{2}{p-1}} (\log |x|)^\frac{r}{p-1}
\end{equation}
as $|x|\to \infty$, give the critical behaviors in the case of \eqref{eq:1.14} and \eqref{eq:1.15}, 
respectively. 
See Section~\ref{section:5} for precise statements.

\medskip 

This paper is organized as follows: 
In Section~\ref{section:2} we show preliminary results. 
In particular, we recall properties of self-similar solutions 
for model equations, that is, power type and exponential type semilinear heat equations. 
In Section~\ref{section:3} and \ref{section:4} we give the proofs of Theorems~\ref{Theorem:1.1}--\ref{Theorem:1.3}. 
Section~\ref{section:5} is devoted to 
applications of our results. 

\section{Preliminaries}
\label{section:2}

In this section we give preliminary results. 
We first study the behavior of the function $f=f(s)$ as $s\to +0$. 
We next recall 
properties of forward self-similar solution for model problems, 
that is, problem~\eqref{eq:1.1} with $f(u) = u^p$ ($p>1$) and $f(u) = e^u$. 
Finally, we introduce the transformations of solutions 
and the equations via a quasi self-similar transformation. 


We first state properties of the functions related to $f$. 
Under the assumptions~\eqref{eq:1.2} and \eqref{eq:1.3}, 
we can easily see that the inverse function $F^{-1}$ of $F$ exists. 
In particular, since $f\in C^1([0,\infty))$ and $f(0) = 0$, 
it follows from \cite[Lemma~3.1]{FI3} that $F(s) \to \infty$ as $s\to +0$. 
This implies that 
\begin{gather}
  \label{eq:2.1}
  (0,\infty)\ni s  \mapsto F(s)\in (0,\infty) 
  \ \text{is a decreasing function},
  \\
  \label{eq:2.2}
  (0,\infty)\ni \sigma  \mapsto F^{-1}(\sigma)\in (0,\infty) 
  \ \text{is a decreasing function}, 
  \\ 
  \label{eq:2.3}
  F^{-1}(\sigma) \to 0 \ \text{as}\ \sigma \to \infty. 
\end{gather} 
We next consider the limit $q$ defined by \eqref{eq:1.11}, 
and prove that it has to satisfy $q\ge 1$. 

\begin{Lemma}
  \label{Lemma:2.1}
  Let $N\ge 1$ and $f\in C^1([0,\infty)) \cap C^2((0,\infty))$ satisfy \eqref{eq:1.2} and \eqref{eq:1.3}.
  Assume that the limit $q$ defined by \eqref{eq:1.11} exists.
  Then $q\ge 1$. 
\end{Lemma}
\noindent
\textbf{Proof}.
Suppose on the contrary that $q<1$.
Then we can take constants $\delta\in (0,1)$ and $\epsilon_0>0$ such that
\begin{equation*}
  f'(s)F(s) \le 1-\delta
\end{equation*}
for all $s\in (0,\epsilon_0]$.
Since $f=-1/F'$, we have $f'=F''/(F')^2$.
These implies that
\begin{equation}
  \notag
  \frac{F''F}{(F')^2} \le 1-\delta
\end{equation}
for all $s\in (0,\epsilon_0]$.
Then, paying attention to the fact that $F'=-1/f<0$, we have
\begin{equation}
  \notag
  \left( \log |F'(s)| \right)'
  = \frac{F''}{F'} \ge (1-\delta)\frac{F'}{F} = (1-\delta) (\log F(s))'.
\end{equation}
Integrating both sides of the above inequality from $s$ to $\epsilon_0$,
we obtain
\begin{equation}
  \notag
  \log \frac{|F'(\epsilon_0)|}{|F'(s)|} \ge (1-\delta) \log \frac{F(\epsilon_0)}{F(s)}
\end{equation}
for all $s\in (0,\epsilon_0]$,
which implies that there exists a constant $C>0$ such that
\begin{equation}
  \notag
  |F'(s)| \le CF(s)^{1-\delta}
  \iff
  F'(s) \ge -CF(s)^{1-\delta}.
\end{equation}
Then we have
\begin{equation}
  \notag
  \left( F(s)^\delta \right)' = \delta F(s)^{\delta-1}F'(s) \ge -C\delta,
\end{equation}
and obtain
\begin{equation}
  \notag
  F(\epsilon_0)^\delta - F(s)^\delta \ge -C(\epsilon_0-s)
  \iff
  F(s)^\delta \le F(\epsilon_0)^\delta + C(\epsilon_0-s).
\end{equation}
This implies that
\begin{equation}
  \notag
  \limsup_{s\to +0} F(s) < \infty,
\end{equation}
but this contradicts the fact that $F(s)\to \infty$ as $s\to +0$ (see \cite[Lemma~3.1]{FI3}).
Therefore, we obtain $q\ge 1$.
\qed 

\medskip

The existence of solutions for problem~\eqref{eq:1.1} 
is proved via a construction of supersolutions (see for example \cite{RS} and \cite{S}).

\begin{Proposition}
  \label{Proposition:2.1}
  Let $f\in C^1([0,\infty)) \cap C^2((0,\infty))$ satisfy \eqref{eq:1.2}.
  Let $u_0\in C(\mathbb{R}^N)$ be a nonnegative function.
  Assume that there exist a supersolution $\overline{u}$ and a subsolution $\underline{u}$
  belonging to $C^{2,1}(\mathbb{R}^N \times (0,T)) \cap C(\mathbb{R}^N \times [0,T))$
  of problem~\eqref{eq:1.1} in $\mathbb{R}^N \times (0,T)$, where $0<T\le \infty$, 
  that is, 
  \begin{equation*}
    \partial_t \overline{u} \ge \Delta \overline{u} + f(\overline{u}) 
    \quad\mbox{in}\,\,\, \mathbb{R}^N \times (0,T), 
    \qquad 
    \overline{u}(x,0) \ge u_0(x) \quad\mbox{in}\,\,\, \mathbb{R}^N, 
  \end{equation*}
  and 
  \begin{equation*}
    \partial_t \underline{u} \le \Delta \underline{u} + f(\underline{u}) 
    \quad\mbox{in}\,\,\, \mathbb{R}^N \times (0,T), 
    \qquad 
    \underline{u}(x,0) \le u_0(x) \quad\mbox{in}\,\,\, \mathbb{R}^N. 
  \end{equation*}
  Then there exists a solution $u\in C^{2,1}(\mathbb{R}^N \times (0,T))\cap C(\mathbb{R}^N \times [0,T))$ of problem~\eqref{eq:1.1}
  satisfying $\underline{u}(x,t) \le u(x,t)\le \overline{u}(x,t)$ in $\mathbb{R}^N \times [0,T)$.
\end{Proposition} 

\begin{Remark}
  \label{Remark:2.1}
  A function $\underline{u}(x,t)\equiv 0$ is a trivial subsolution of problem~\eqref{eq:1.1}. 
  Once a supersolution $\overline{u}$ is constructed, 
  there exists a solution $u$ if problem~\eqref{eq:1.1}  
  satisfying $0 \le u(x,t)\le \overline{u}(x,t)$. 
\end{Remark}

We next study a sufficient condition for the global in time existence of solutions of problem~\eqref{eq:1.1}. 
In order to prove the nonexistence of global solutions,
we use the following lemma.
Put 
\begin{equation}
  \label{eq:2.4}
  (e^{t\Delta}u)(x) 
  = (4\pi t)^{-\frac{N}{2}} \int_{\mathbb{R}^N} e^{-\frac{|x-y|^2}{4t}} u_0(y) \, dy 
\end{equation}
in $(x,t) \in \mathbb{R}^N \times (0,\infty)$ 
for $u_0\in L^\infty (\mathbb{R}^N)$. 

\begin{Lemma}
  \label{Lemma:2.2}
  Let $N\ge 1$ and $f\in C^1([0,\infty)) \cap C^2((0,\infty))$ satisfy \eqref{eq:1.2} and \eqref{eq:1.3}.
  Let $u_0\in L^\infty(\mathbb{R}^N)$ be a nonnegative initial function.
  If there exists a global in time solution of problem~\eqref{eq:1.1},
  then
  \begin{equation}
    \notag
    F(\|e^{t\Delta}u_0\|_{L^\infty(\mathbb{R}^N)}) \ge t
  \end{equation}
  for all $t>0$.
\end{Lemma}

\noindent
\textbf{Proof}.
We may assume without loss of generality that $u_0\not\equiv 0$.
Since the function
\begin{equation*}
  \underline{u}(x,t) := F^{-1}
  \left[
  F((e^{t\Delta}u_0)(x)) - t
  \right]
\end{equation*}
is a subsolution of problem~\eqref{eq:1.1}
(see for example the proof of \cite[Theorem~1.2]{FI3}),
we have $\underline{u}(x,t)\le u(x,t)$
and $\underline{u}(x,t) < \infty$ for all $(x,t)\in \mathbb{R}^N\times (0,\infty)$
if the solution exists globally in time.
Then, since it follows from \eqref{eq:1.3} that
$F(u)\to 0$ as $u\to\infty$, that is,
$F^{-1}(u)\to \infty$ as $u\to +0$,
we see that $F((e^{t\Delta}u_0)(x)) > t$,
which yields the desired inequality.
$\qed$

\medskip 

We recall forward self-similar solutions
for problem~\eqref{eq:1.1} with $f(u) = u^p$, where $p>1$.
In particular, we introduce the result on the solution set of forward self-similar solutions
which has been obtained by Naito~\cite{N}. 
Consider a power type semilinear heat equation
\begin{equation}
  \label{eq:2.5}
  \partial_t u = \Delta u + u^p,
  \qquad x\in \mathbb{R}^N, \,\,\, t>0,
\end{equation}
where $p>1$.
Forward self-similar solution is a solution $u$ of \eqref{eq:2.5}
with the form
\begin{equation}
  \notag
  u(x,t) = (t+1)^{-\frac{1}{p-1}}v(y),
  \qquad
  y = \frac{x}{\sqrt{t+1}}.
\end{equation}
If $u$ is a forward self-similar solution of \eqref{eq:2.5},
then $v$ has to satisfy
\begin{equation}
  \label{eq:2.6}
  \Delta v + \frac{1}{2}y\cdot \nabla v + \frac{1}{p-1}v + v^p = 0
  \quad\mbox{in}\,\,\, \mathbb{R}^N.
\end{equation}
Now we focus on radially symmetric solutions of \eqref{eq:2.6} with $v(0)=\alpha>0$
and $|y|^{2/(p-1)}v(y) \to \ell>0$ as $|y|\to\infty$.
Let $v=v(r)$ and $r=|y|$.
Then $v$ satisfies $v(0) = \alpha$ and
\begin{equation}
  \label{eq:2.7}
  \begin{cases}
    v'' + \dfrac{N-1}{r}v' + \dfrac{1}{2}rv' + \dfrac{1}{p-1}v + v^p = 0,
    & r>0,
    \\[5pt]
    v'(0) = 0,
    \quad
    v(0) = \alpha.
  \end{cases}
\end{equation}
We denote by $v_\alpha$ the unique solution of \eqref{eq:2.7}.
Furthermore, the limit
\begin{equation}
  \notag
  \ell_p(\alpha) := \lim_{r\to\infty} r^\frac{2}{p-1}v_\alpha(r)
  \in [0,\infty)
\end{equation}
exists.
The properties of a solution for \eqref{eq:2.7} and $\ell_p(\alpha)$ have been studied
in \cite{N},
and we have the following:
\begin{Proposition}
  \label{Proposition:2.2}
  Let $N\ge 1$ and $p>1+2/N$.
  There exists a constant $\alpha_*\in (0,\infty]$ such that
  \begin{enumerate}
    \item[{\rm (i)}]
    if $\alpha \in (0,\alpha_*)$,
    then $v_\alpha(r)>0$ in $[0,\infty)$ and $\ell_p(\alpha)>0$.

    \item[{\rm (ii)}]
    $\ell(\alpha)$ is strictly increasing in $\alpha\in (0,\alpha_*)$
    and $\ell_p(\alpha)\to 0$ as $\alpha\to +0$.
  \end{enumerate}
\end{Proposition}

We see from \cite[Proposition~3.4]{HW} that
a universal decay estimate for solutions of \eqref{eq:2.7} follows.

\begin{Proposition}
  \label{Proposition:2.3}
  Let $N\ge 1$ and $p>1+2/N$.
  Let $\alpha_*$ be the constant given in Proposition~{\rm \ref{Proposition:2.2}}
  and $\alpha_0<\alpha_*$.
  Then there exists a positive constant $C$, depending only on $N$, $p$ and $\alpha_0$, such that
  \begin{equation*}
    \sup_{r\ge 0}\ (1+r)^\frac{2}{p-1} v_\alpha (r)
    \le
    C\alpha
  \end{equation*}
  for all $\alpha\in (0,\alpha_0)$.
\end{Proposition}

We also recall properties of forward self-similar solutions for problem~\eqref{eq:1.1}
with $f(u) = e^u$.
Consider 
\begin{equation}
  \label{eq:2.8}
  \partial_t u = \Delta u + e^u, 
  \qquad x\in \mathbb{R}^N, \,\,\, t>0. 
\end{equation}
A function $u$ is said to be a forward self-similar solution of \eqref{eq:2.8} if 
there exists a function $v$ such that 
\begin{equation*}
  u(x,t) = \log \frac{1}{ t+1} + v(y), 
  \qquad 
  y = \frac{x}{\sqrt{t+1}}, 
\end{equation*}
where $v$ satisfies 
\begin{equation*}
  \Delta v + \frac{1}{2}y\cdot \nabla v + 1 + e^v = 0 
  \quad\mbox{in}\quad \mathbb{R}^N. 
\end{equation*}
As in the case of a power type semilinear heat equation, 
we are interested in radially symmetric forward self-similar solutions, 
so we consider 
\begin{equation}
  \label{eq:2.9}
  \begin{cases}
    v''+\dfrac{N-1}{r}v'+\dfrac{1}{2}rv'+1+e^v=0,
    & r>0,
    \\[5pt]
    v'(0) = 0,
    \quad
    v(0) = \alpha, 
  \end{cases}
\end{equation}
satisfying 
\begin{equation}
  \notag
  \ell_{\operatorname{exp}}(\alpha)
  :=
  \lim_{r\to \infty}
  ( 2\log r + v_\alpha(r) )
  \in \mathbb{R}.
\end{equation}
Denote the unique solution of \eqref{eq:2.9} by $v_\alpha$. 
We have the following proposition (see \cite{F2}). 

\begin{Proposition}
  \label{Proposition:2.4}
  Let $N\ge 1$ and $\alpha_0\in \mathbb{R}$.
  There exists a constant $C$, depending only on $N$ and $\alpha_0$, such that
  \begin{equation}
    \notag
    \sup_{r\ge 0} \Big( 2\log (1+r) + v_\alpha(r) \Big)
    \le C + \alpha
  \end{equation}
  for all $\alpha<\alpha_0$.
\end{Proposition}

\begin{Remark}
  \label{Remark:2.2}
  If there exists a positive solution $v_\alpha$ of \eqref{eq:2.7},
  then it follows that
  \begin{equation}
    \notag
    \left( r^{N-1}e^{\frac{r^2}{4}}v_\alpha'(r) \right)'
    =
    r^{N-1}e^{\frac{r^2}{4}} \left( v_\alpha'' + \frac{N-1}{r}v_\alpha' + \frac{r}{2}v_\alpha' \right)
    =
    -r^{N-1}e^{\frac{r^2}{4}} \left( \frac{1}{p-1}v_\alpha + v_\alpha^p \right)
    \le 0
  \end{equation}
  for $r>0$. 
  Thus we have $r^{N-1}e^{\frac{r^2}{4}}v_\alpha'(r) \le r^{N-1}e^{\frac{r^2}{4}}v_\alpha'(r)\Big|_{r=0} = 0$,
  that is, $v'(r)\le 0$ for $r>0$.
  In particular, we obtain $v_\alpha(r)\le v_\alpha(0)=\alpha$. 
  Similarly, 
    if there exists a solution $v_\alpha$ of \eqref{eq:2.9},
  we have
  \[
      \left( r^{N-1}e^{\frac{r^2}{4}}v_\alpha'(r) \right)'
      =-r^{N-1}e^{\frac{r^2}{4}}(1+e^v)\le 0
  \]
  and hence 
  $v_\alpha(r)\le v_\alpha(0)=\alpha$.
\end{Remark}

We finally introduce key transformations based on \eqref{eq:1.7}. 
The proofs of main theorems rely on the following proposition:

\begin{Proposition}
  \label{Proposition:2.5}
  Let $f$, $g\in C^1([0,\infty)) \cap C^2((0,\infty))$ be functions satisfying \eqref{eq:1.2} and \eqref{eq:1.3}, and set
  \begin{equation*}
    G(s):=\int_s^{\infty}\frac{d\eta}{g(\eta)}
  \end{equation*}
  for $s>0$. 
  Define 
  \begin{equation}
    \label{eq:2.10}
    u(x,t)
    :=
    F^{-1}\left[(t+1)F\left(w\left(\frac{|x|}{\sqrt{t+1}}\right)\right)\right]
    \quad \text{and} \quad 
    w(r):=
    F^{-1}\left[ 
      G(v(r))
    \right]
  \end{equation}
  for $v\in C^{2}([0,\infty))$.
  Then $u$, $v$ and $w$ satisfy the following 
  equalities:
  \begin{equation}
    \begin{aligned}
      \label{eq:2.11}
      &
      \frac{1}{f(u)}
      \left[
        \partial_t u -\Delta u -f(u)+\frac{f'(u)}{f(u)}|\nabla u|^2
      \right]
      \\
      &
      =
      \frac{1}{f(w)}
      \left[
        -w''-\frac{N-1}{r}w-\frac{r}{2}w' -f(w)-f(w)F(w)+\frac{f'(w)}{f(w)}(w')^2
      \right]
      \\
      &
      =
      \frac{1}{g(v)}
      \left[
        -v''-\frac{N-1}{r}v-\frac{r}{2}v' -g(v)-g(v)G(v)+\frac{g'(v)}{g(v)}(v')^2
      \right]
    \end{aligned}
  \end{equation}
  for $x\in \mathbb{R}^N$ and $t>0$. 
\end{Proposition}

As a corollary of Proposition~\ref{Proposition:2.5}, 
we have the following result. 

\begin{Corollary}
  \label{Corollary:2.1}
  Assume the same conditions as in Proposition~{\rm \ref{Proposition:2.5}}. 
  If $v\in C^2([0,\infty))$ satisfies 
  \begin{equation}
    \notag 
    v''+\frac{N-1}{r}v+\frac{r}{2}v' +g(v)+g(v)G(v)=0,
  \end{equation}
  then 
  \[
    u(x,t)=F^{-1}\left[(t+1)G\left(\frac{|x|}{\sqrt{t+1}}\right)\right]
  \]
  satisfies
  \begin{equation}
    \notag 
    \partial_t u -\Delta u=f(u)+\frac{|\nabla u|^2}{f(u)F(u)}
    \Bigl[ g'\bigl(v(r)\bigr)G\bigl(v(r)\bigr)-f'\bigl(u(x,t)\bigr)F\bigl(u(x,t)\bigr)\Bigr],
    \quad
    r=\frac{|x|}{\sqrt{t+1}}, 
  \end{equation}
  for all $x\in \mathbb{R}^N$ and $t>0$.
\end{Corollary}

Corollary~\ref{Corollary:2.1} follows from \eqref{eq:2.11} and the facts that
\[
  |\nabla u(x,t)|=\sqrt{t+1}\cdot \frac{f\bigl(u(x,t)\bigr)}{g\bigl(w(r)\bigr)}
  |w'(r)|
  \quad \text{and} \quad 
  t+1=\frac{F\bigl(u(x,t)\bigr)}{G\bigl(w(r)\bigr)}
  \quad \text{for} \quad 
  r=\frac{|x|}{\sqrt{t+1}}.
\]
We apply Corollary~\ref{Corollary:2.1} with $g(v)=v^p$ or $g(v)=e^v$ 
for the proofs of main theorems. 
The rest of this section is devoted to the proof of Proposition~\ref{Proposition:2.5}. 

\medskip

\noindent
\textbf{Proof of Proposition~\ref{Proposition:2.5}}.
Let $v\in C^2([0,\infty))$ and 
$u$, $w$ be the functions defined by \eqref{eq:2.10}.
Then $u\in C^{2,1}(\mathbb{R}^N\times (0,\infty))$ and $w\in C^2([0,\infty))$. 
Let $y:=\frac{x}{\sqrt{t+1}}$ and $r=|y|$.
For the sake of simplicity, we use the same notation $v(x):=v(|x|)$ and $w(x):=w(|x|)$
for $x\in \mathbb{R}^N$.
We prove the first equality.
It follows from direct computation that
\[
  \begin{aligned}
    \partial_t u(x,t) 
    &
    =-f(u)\left[F(w)+\frac{1}{f(w)}\frac{y}{2}\cdot \nabla w(y)\right],
    \\
    \nabla u(x,t)
    &
    =\frac{f(u)}{f(w)}\sqrt{t+1}\nabla w(y),
    \\
    \Delta u(x,t)
    &
    =\frac{f(u)}{f(w)}\Delta w +\frac{f'(u)f(w)\sqrt{t+1}\nabla u \cdot \nabla w-f(u)f'(v)|\nabla w(y)|^2}{f(w)^2}
    \\
    &
    =\frac{f(u)}{f(w)}\Delta w 
    +\frac{f'(u)}{f(u)}|\nabla u|^2
    -\frac{f'(w)f(u)}{f(w)^2}|\nabla w|^2.
  \end{aligned}
\]
Therefore, it holds
\[
  \frac{1}{f(u)}\left(\partial_t u-\Delta u+\frac{f'(u)}{f(u)^2}|\nabla u|^2 \right)
  =
  -F(w)-\frac{1}{f(w)}\frac{y}{2}\cdot \nabla w
  -\frac{\Delta w}{f(w)}+ \frac{f'(w)}{f(w)^2}|\nabla w|^2.
\]
This yields the first equality.
We can prove the second equality in the same manner. 
It holds
\begin{equation}
  \notag 
  \begin{aligned}
    \nabla w
    &
    =\frac{f(w)}{g(v)}\nabla v,
    \\
    \Delta w 
    &
    =\frac{f(w)}{g(v)}\Delta v+\frac{f'(w)g(v)\nabla w\cdot \nabla v-f(w)g'(v)|\nabla v|^2}{g(v)^2}
    \\
    &
    =\frac{f(w)}{g(v)}\Delta v
    +\frac{f'(w)}{f(w)}|\nabla w|^2-\frac{f(w)g'(v)}{g(v)^2}|\nabla v|^2
    \frac{f'(w)g(v)\nabla w\cdot \nabla v-f(w)g'(v)|\nabla v|^2}{g(v)^2}.
  \end{aligned}
\end{equation}
This together with 
\begin{equation*}
  \dfrac{y\cdot \nabla w}{f(w)}=\dfrac{y\cdot \nabla v}{g(v)}, \qquad F(w)=G(v), 
\end{equation*}
gives the desired equality. 
Thus we complete the proof of Proposition~\ref{Proposition:2.5}.
\qed

\section{Proofs of Theorems~\ref{Theorem:1.1} and \ref{Theorem:1.2}}
\label{section:3}

This section is devoted to the proofs of Theorems~\ref{Theorem:1.1} and \ref{Theorem:1.2}.

\medskip

\noindent
\textbf{Proof of Theorem~\ref{Theorem:1.1}}.
Under the assumptions of Theorem~\ref{Theorem:1.1}, 
by \eqref{eq:1.12} we can take a constant $q^*\in (q,1+N/2)$.
Observe that $q^*>1$.
Then there exists a constant $s^*>0$ such that
\begin{equation}
  \label{eq:3.1}
  f'(s)F(s) \le q^*
\end{equation}
for all $s\in (0,s^*)$.

In order to construct a global in time solution of problem~\eqref{eq:1.1},
we apply the transformation in
Corollary~\ref{Corollary:2.1}
with
$g(s)=s^{\frac{q^*}{q^*-1}}$.
Let $p^*:=\frac{q^*}{q^*-1}$.
Then it holds that
\begin{equation}
  \label{eq:3.2}
  g(s)=s^{p^*},\ \ 
  G(s)=\frac{1}{p^*-1}s^{-(p^*-1)},\ \ g'(s)G(s)=q^*=\frac{p^*}{p^*-1}, \ \ 
  g(s)G(s)=\frac{1}{p^*-1}s, 
\end{equation}
for all $s>0$.
Let $\alpha$ be a positive constant to be chosen later and $v_\alpha$ be the solution of
\begin{equation}
  \label{eq:3.3}
  \begin{cases}
    v'' + \dfrac{N-1}{r}v' + \dfrac{1}{2}rv' + \dfrac{1}{p^*-1}v + v^{p^*} = 0,
    & r>0,
    \\[5pt]
    v'(0) = 0,
    \quad
    v(0) = \alpha. 
  \end{cases}
\end{equation}
It follows from Proposition~\ref{Proposition:2.2} that
\begin{equation}
  \label{eq:3.4}
  |x|^{\frac{2}{p^*-1}}v_{\alpha}(|x|)
  \to c_\alpha
\end{equation}
as $|x|\to \infty$ for some constant $c_\alpha\in (0,\infty)$.
Define
\[
  u_{\alpha}(x,t):=F^{-1}\left[
    (t+1)\cdot \frac{1}{p^*-1}v_{\alpha}\left( 
      \frac{|x|}{\sqrt{t+1}}
    \right)^{-(p^*-1)}
  \right] 
  \quad\mbox{in}\,\,\, 
  (x,t) \in \mathbb{R}^N \times [0,\infty). 
\] 
Then, by Corollary~\ref{Corollary:2.1} with \eqref{eq:3.2} and \eqref{eq:3.3} we see that $u_{\alpha}$
satisfies
\[
  \left\{
    \begin{aligned}
      \partial_t u_{\alpha}-\Delta u_{\alpha}
      &
      =f(u_{\alpha})+\frac{|\nabla u_{\alpha}|^2}{f(u_{\alpha})F(u_{\alpha})}\Bigl(q^*-f'(u_{\alpha})F(u_{\alpha})\Bigr),
      \\
      u_{\alpha}(x,0)
      &
      =
      F^{-1}\left[\frac{1}{p^*-1}v_{\alpha}(|x|)^{-(p^*-1)}\right]. 
    \end{aligned}
  \right.
\]
Since $v_{\alpha}(r)\le \alpha$ for all $r>0$ (see Remark~\ref{Remark:2.2}), 
by \eqref{eq:2.2} we have 
\[
  u_{\alpha}(x,t)\le F^{-1}\left[\frac{1}{p^*-1}\alpha^{-(p^*-1)}\right],
\]
hence $u_{\alpha}$ is bounded. 
Furthermore, taking a small enough $\alpha>0$, 
by \eqref{eq:2.3} 
we obtain $u_{\alpha}(x,t)\le s^*$ in $\mathbb{R}^N \times [0,\infty)$.
On the other hand, 
by \eqref{eq:2.2} and \eqref{eq:3.4} we can take a constant $\gamma_*>0$ such that 
\begin{equation}
  \label{eq:3.5}
  \frac{v_\alpha^{-(p^*-1)} }{p^*-1}
  \le 
  \gamma^{-1}(|x|^2+1)
  \iff 
  F^{-1}\left[\frac{1}{p^*-1}v_{\alpha}(|x|)^{-(p^*-1)}\right]
  \ge 
  F^{-1} \Bigl[
    \gamma^{-1}(|x|^2+1)
  \Bigr]
  =u_0
\end{equation}
in $\mathbb{R}^N$, where $\gamma\in (0,\gamma_*)$. 
Therefore, we see from \eqref{eq:3.1} that 
the function $u_{\alpha}$ is a global in time supersolution of 
\begin{equation}
  \notag 
  \left\{
    \begin{aligned}
      \partial_t u-\Delta u
      &
      =f(u), && x\in \mathbb{R}^N,\ t>0,
      \\
      u(x,0)
      &
      =
      F^{-1}\left[\frac{1}{p^*-1}v_{\alpha}(|x|)^{-(p^*-1)}\right], 
      && x\in \mathbb{R}^N.
    \end{aligned}
  \right.
\end{equation}
Thanks to Proposition~\ref{Proposition:2.1} and \eqref{eq:3.5}, 
one can construct a global solution $u$ of problem~\eqref{eq:1.1} satisfying
\begin{equation}
  \label{eq:3.6}
  u(x,t)\le u_{\alpha}(x,t)
  =
  F^{-1}\left[(t+1) \cdot \frac{1}{p^*-1}v_{\alpha}\left(\frac{|x|}{\sqrt{t+1}}\right)^{-(p^*-1)}\right].
\end{equation}
This completes the proof of Theorem~\ref{Theorem:1.1}.
\qed

\medskip

\noindent
\textbf{Proof of Theorem~\ref{Theorem:1.2}}.
We use the same notation in the proof of Theorem~\ref{Theorem:1.1}.
Since $F$ is decreasing, by \eqref{eq:3.6} we get an upper bound of $u$
and obtain
\begin{equation}
  \notag 
  \frac{|x|^2}{F(u(x,t))}
  \le \frac{|x|^2}{t+1} 
  (p^*-1)v_{\alpha}\left(\frac{|x|}{\sqrt{t+1}}\right)^{p^*-1}
  =\frac{r^2}
  {\frac{1}{p^*-1}v_{\alpha}(r)^{-(p^*-1)}},
  \quad r=\frac{|x|}{\sqrt{t+1}}.
\end{equation}
Taking
$W^*$ as 
\[
  W^*(y):=F^{-1}\left[\frac{1}{p^*-1}v_{\alpha}(|y|)^{-(p^*-1)}\right],
\]
we obtain the upper estimate of \eqref{eq:1.13}.
Furthermore, 
by putting $g(s):=s^{p^*}$ and ${G}(s)=\int_s^{\infty}\frac{1}{g(\eta)}d\eta=\frac{1}{p^*-1}s^{-(p^*-1)}$, we have
$g'(s)G(s)=\frac{p^*}{p^*-1}=q^*$ and
$W^*=F^{-1}\left[G(v_{\alpha})\right]$.
Therefore Proposition~\ref{Proposition:2.5} and \eqref{eq:3.3} gives us that
$W^*$ satisfies
\[
  \Delta W^*+\frac{y}{2}\cdot \nabla W^*(y) 
  +f(W^*)F(W^*)+f(W^*)
  +\frac{|\nabla W^*|^2}{f(W^*)F(W^*)}\Bigl[q^*-f'(W^*)F(W^*)\Bigr]=0 
\]
in $\mathbb{R}^N$. 

In the rest of the proof, we obtain a lower bound of the solution.
Let $\beta>0$ and 
\begin{equation*}
  \beta_q = 
  \begin{cases}
    \beta & \mbox{if} \ \ q>1, 
      \\
      \log \beta & \mbox{if} \ \ q=1. 
  \end{cases}
\end{equation*}
Observe that $\beta_1 \to -\infty$ as $\beta\to +0$. 
If $q>1$, there exist $s_*>0$ and $q_*\in (1,q)$ such that 
$f'(s)F(s)\ge q_*$ for all $s\in (0,s_*]$. 
On the other hand, for the case $q=1$ we assume that there exists $s_*>0$ 
such that $f'(s)F(s)\ge 1$ for all $s\in (0,s_*]$.
For the lower estimate of \eqref{eq:1.13}, we consider
\begin{equation}
  \label{eq:3.7}
  g(s)
  :=
  \left\{
    \begin{aligned}
      & 
      s^{p_*}&& \text{for}\ s>0 && \text{if}\ \ q>1,
      \\
      & e^s && \text{for}\ s\in \mathbb{R}&& \text{if}\ \ q=1,
    \end{aligned}
  \right.
\end{equation}
where $p_*:=\frac{q_*}{q_*-1}$. 
It is easy to check that
if $q>1$ then
\[
  G(s)
  =
  \frac{1}{p_*-1}s^{-(p_*-1)}
  \quad 
  \text{and}
  \quad 
  g(s)G(s)=\frac{1}{p_*-1}s
  \quad \text{for\ all\ }s>0, 
\]
and if $q=1$ then
\[
  G(s)
  =
  e^{-s}
  \quad 
  \text{and}
  \quad 
  g(s)G(s)=1
  \quad \text{for\ all\ }s\in \mathbb{R}.
\]
Hence, the problem 
\begin{equation}
  \label{eq:3.8}
  \left\{
    \begin{aligned}
      &
      v''+\frac{N-1}{r}v'+\frac{1}{2}rv'+g(v)G(v)+g(v)=0,
      \quad r>0, 
      \\
      &
      v'(0)=0,\quad v(0)=\beta_q,
    \end{aligned}
  \right.
\end{equation}
coincides with the model problem 
\eqref{eq:2.7} with $\alpha=\beta_q$
and 
\eqref{eq:2.9} with $\alpha=\beta_1$
for $q>1$ and $q=1$, respectively. 
Remark that the solution of \eqref{eq:3.8} uniquely exists with $g$ defined by \eqref{eq:3.7}.
Let $v_{\beta}$ be the solution of \eqref{eq:3.8}.
Define
\[
  u_{\beta}(x,t)
  :=
  F^{-1}
  \left[
   (t+1)G(v_{\beta}(r))
  \right],
  \qquad
  r=\frac{|x|}{\sqrt{t+1}}. 
\]
Then, by Corollary~\ref{Corollary:2.1} 
we see that $u_\beta$ satisfies
\begin{equation}
  \label{eq:3.9}
  \left\{
    \begin{aligned}
      &
      \partial_t u_{\beta} -\Delta u_{\beta}=f(u_{\beta})+\frac{|\nabla u_{\beta}|^2}{f(u_{\beta})F(u_{\beta})}
      \Bigl(g'\bigl(v_{\beta}(r)\bigr)G\bigl(v_{\beta}(r)\bigr)-f'\bigl(u_{\beta}(x,t)\bigr)F\bigl(u_{\beta}(x,t)\bigr)\Bigr),
      \\
      &
      u_{\beta}(x,0)=F^{-1}[G(v_{\beta}(|x|))],
    \end{aligned}
  \right.
\end{equation}
in $\mathbb{R}^N\times [0,\infty)$.
We show that $u_{\beta}$ is a subsolution of problem~\eqref{eq:1.1}.
Recall that
\begin{equation}
  \label{eq:3.10}
  g'(s)G(s)
  =
  \left\{
    \begin{aligned}
      & q_* && \text{if}\ \ q>1,
      \\
      &1 && \text{if}\ \ q=1.
    \end{aligned}
  \right.
\end{equation}
Since it follows from Remark~\ref{Remark:2.2} that $v_{\beta}\le \beta$, 
by \eqref{eq:2.1} and \eqref{eq:2.2} we have
\[
  u_{\beta}(x,t)\le 
  F^{-1}
  \left[
   (t+1)G(v_{\beta}(r))
  \right]\le F^{-1}[G(\beta_q)].
\]
This together with \eqref{eq:2.3} implies that there exists $\beta_*>0$ such that
$u_{\beta}(x,t)\le s_*$ in $\mathbb{R}^N \times [0,\infty)$ for $\beta\in (0,\beta_*)$ 
and hence
\[
  f'(u_{\beta})F(u_{\beta})\ge 
  \left\{
    \begin{aligned}
      & q_* && \text{if}\ \ q>1,
      \\
      &1 && \text{if}\ \ q=1.
    \end{aligned}
  \right.
\]
in $\mathbb{R}^N\times [0,\infty)$.
Therefore, by \eqref{eq:3.9} and \eqref{eq:3.10} we have 
\begin{equation}
  \label{eq:3.11}
  \partial_t u_{\beta}-\Delta u_{\beta}\le f(u_{\beta})
\end{equation}
in $\mathbb{R}^N \times (0,\infty)$. 

Next we focus on the initial data $u_{\beta}(x,0)=F^{-1}[G(v_{\beta}(|x|))]$. 
By Propositions~\ref{Proposition:2.3} and \ref{Proposition:2.4}, 
$v_{\beta}$ satisfies
\begin{equation}
  \label{eq:3.12}
  v_{\beta}(r)
  \le
  \left\{
    \begin{aligned}
      &C_1\beta (1+r)^{-\frac{2}{p_*-1}} && \text{if}\ \ q>1,
      \\
      &C_2+\beta_1 -2\log(1+r) && \text{if}\ \ q=1, 
    \end{aligned}
  \right.
\end{equation}
where $C_1>0$ and $C_2\in \mathbb{R}$ are constants. 
Let $\gamma_*>0$ be the constant given in Theorem~\ref{Theorem:1.1} 
and fix $\gamma\in (0,\gamma_*)$. 
Then, taking a small enough $\beta_*>0$ if necessary, 
by \eqref{eq:3.12} we obtain 
\[
  G(v_{\beta})\ge 
  \left\{
    \begin{aligned}
      &
      (C_1\beta)^{-(p_*-1)}(1+r)^2 
      \ge 
      \gamma^{-1}(r+1)^2
      && \text{if}\ \ q>1,
      \\
      &
      e^{-(C_2+\beta_1)}(1+r)^2 
      \ge 
      \gamma^{-1}(r+1)^2
      && \text{if}\ \ q=1, 
    \end{aligned}
  \right.
\]
for $\beta\in (0,\beta_*)$. 
Observe again that $\beta_1\to -\infty$ as $\beta\to +0$. 
Therefore, it holds that
\[
  u_{\beta}(x,0)
  =
  F^{-1}\bigl[G(v_{\beta}(|x|))\bigr]
  \le
  F^{-1}\left[\gamma^{-1}(|x|^2+1)\right]=u_0(x). 
\]
This together with \eqref{eq:3.11} implies that $u_{\beta}$ is a subsolution of problem~\eqref{eq:1.1}.
Then we have 
\[
  u(x,t)
  \ge 
  u_{\beta}(x,t)
  =
  F^{-1}
  \left[
    (t+1)G(v_{\beta}(r))
  \right],
  \quad
  r=\frac{|x|}{\sqrt{t+1}},
\]
and equivalently 
\[
  \frac{|x|^2}{F\bigl(u(x,t)\bigr)}\ge \frac{|y|^2}{G\bigl(v_{\beta}(|y|)\bigr)},
  \qquad 
  y=\frac{x}{\sqrt{t+1}},
\]
hence we obtain the lower estimate of \eqref{eq:1.13}
by choosing $W_*$ as 
\[
  W_*(y):=F^{-1}\Bigl[G\bigl(v_{\beta}(|y|)\bigr)\Bigr].
\]
Moreover, in view of Proposition~{\rm \ref{Proposition:2.5}},
the function $W_*$ satisfies
\[
  \Delta W_*+\frac{y}{2}\cdot \nabla W_*(y)
  +f(W_*)F(W_*)+f(W_*)
  +\frac{|\nabla W_*|^2}{f(W_*)F(W_*)}\Bigl[q_*-f'(W_*)F(W_*)\Bigr]=0
\]
if $q>1$. 
In the case $q=1$, $W_*$ satisfies the same equation replacing $q_*$ by 1. 
This completes the proof of Theorem~\ref{Theorem:1.2}.
\qed

\section{Proof of Theorem~\ref{Theorem:1.3}}
\label{section:4}

In this section we prove Theorem~\ref{Theorem:1.3}.
%

\medskip

\noindent
\textbf{Proof of Theorem~\ref{Theorem:1.3}
}.
Since all nontrivial nonnegative solutions blow up in finite time if $q>1+N/2$,
we assume that $1\le q\le 1+N/2$.
Then we can take a constant $L>0$ such that
\begin{equation}
  \label{eq:4.1}
  \frac{|x|^2}{F(u_0(x))} \ge \gamma
  \iff
  u_0(x) \ge F^{-1}(\gamma^{-1}|x|^2)
\end{equation}
for all $|x|\ge L$.
We prove that if $\gamma$ is large enough, then a global in time solution cannot exist.
Suppose on the contrary that a solution exists globally in time.

If $q>1$, there exist $s_*>0$ and $q_*\in (1,q)$ such that 
$f'(s)F(s)\ge q_*$ for all $s\in (0,s_*]$. 
In the case of $q=1$, we further assume that
there exists $s_*>0$ such that
$f'(s)F(s)\ge 1$ for all $s\in (0,s_*]$.
Let $g$ be the function defined by \eqref{eq:3.7}.
Define $\varphi(s):=F^{-1}\bigl[G(s)\bigr]$. 
It follows from $F(\varphi(s))=G(s)$ that
\[
  \begin{aligned}
    \varphi'(s)
    &
    =\frac{f(\varphi(s))}{g(s)}>0,
    \\
    \varphi''(s)
    &
    =
    \left(\frac{f(\varphi(s))}{g(s)}\right)'
    =
    \frac{f(\varphi(s))}{g(s)^2}\Bigl[f'(\varphi(s))-g'(s)\Bigr]
    \\
    &
    =
    \frac{f(\varphi(s))}{g(s)^2F(\varphi(s))}\Bigl[f'(\varphi(s))F(\varphi(s))-g'(s)G(s)\Bigr]\ge 0, 
  \end{aligned}
\]
for all $s\in (0,s_*]$, 
which implies that $(0,s_*] \ni s \mapsto \varphi(s)$ is convex. 
Without loss of generality, we may assume that $u_0\in L^1(\mathbb{R}^N)\cap L^\infty(\mathbb{R}^N)$ 
and $\|u_0\|_{L^\infty(\mathbb{R}^N)}$ is small enough,
so $0\le u_0(x) \le s_*$ for all $x\in \mathbb{R}^N$.
Thus we see from the Jensen inequality that
\begin{equation}
  \notag 
  (e^{t\Delta}u_0)(x) \ge \varphi\left( (e^{t\Delta}\varphi^{-1}(u_0))(x) \right).
\end{equation}
Then, since $F$ is decreasing, we have 
\begin{equation}
  \label{eq:4.2}
  F\Bigl[
    \|e^{t\Delta}u_0\|_{L^{\infty}(\mathbb{R}^N)}
  \Bigr]
  \le
  F
  \left[
    \varphi
    \Bigl(
      \|e^{t\Delta}\varphi^{-1}(u_0)\|_{L^{\infty}(\mathbb{R}^N)}
    \Bigr)
  \right]
  =G\Bigl(
    \|e^{t\Delta}\varphi^{-1}(u_0)\|_{L^{\infty}(\mathbb{R}^N)}
  \Bigr).
\end{equation}
It remains to estimate $\|e^{t\Delta}\varphi^{-1}(u_0)\|_{L^{\infty}(\mathbb{R}^N)}$ from below.
Since $\varphi^{-1}(s) = G^{-1}\bigl[F(s)\bigr]$, 
by \eqref{eq:2.4} and \eqref{eq:4.1} we have 
\begin{equation}
  \label{eq:4.3}
  \begin{aligned}
    \|e^{t\Delta}\varphi^{-1}(u_0)\|_{L^{\infty}(\mathbb{R}^N)}
    &
    \ge 
    (4\pi t)^{-\frac{N}{2}}\int_{\mathbb{R}^N}{e^{-\frac{|y|^2}{4t}}}\varphi^{-1}(u_0(y))dy
    \\
    &
    \ge
    (4\pi t)^{-\frac{N}{2}}\int_{|y|\ge L}{e^{-\frac{|y|^2}{4t}}}G^{-1}(\gamma^{-1}|y|^2)dy
    \\
    &
    \ge
    (4\pi t)^{-\frac{N}{2}}\int_{\mathbb{R}^N}{e^{-\frac{|y|^2}{4t}}}G^{-1}(\gamma^{-1}|y|^2)dy
    +O(t^{-N/2}).
  \end{aligned}
\end{equation}
If $q>1$, since $G^{-1}(s)=(p_*-1)^{-\frac{1}{p_*-1}}s^{-\frac{1}{p_*-1}}$, we have
\[
  \begin{aligned}
    (4\pi t)^{-\frac{N}{2}}\int_{\mathbb{R}^N}{e^{-\frac{|y|^2}{4t}}}G^{-1}(\gamma^{-1}|y|^2) \, dy
    &
    \ge 
    (4\pi t)^{-\frac{N}{2}}\int_{\mathbb{R}^N}{e^{-\frac{|y|^2}{4t}}}
    (p_*-1)^{-\frac{1}{p_*-1}}
    \gamma^{\frac{1}{p_*-1}}
    |y|^{-\frac{2}{p_*-1}} \, dy
    \\
    &
    = 
    c_1 \gamma^{\frac{1}{p_*-1}}
    (p_*-1)^{-\frac{1}{p_*-1}}t^{-\frac{1}{p_*-1}}
    = 
    c_1 2^{-\frac{1}{p_*-1}} \gamma^{\frac{1}{p_*-1}} G^{-1}\left(\frac{t}{2}\right)
  \end{aligned}
\]
for $t>0$, where $c_1>0$ is a constant.
If $q=1$, since $G^{-1}(s)=\log \frac{1}{s}$, we have
\[
  \begin{aligned}
    (4\pi t)^{-\frac{N}{2}}\int_{\mathbb{R}^N}{e^{-\frac{|y|^2}{4t}}}G^{-1}(\gamma^{-1}|y|^2) \, dy
    &
    \ge 
    (4\pi t)^{-\frac{N}{2}}\int_{\mathbb{R}^N}{e^{-\frac{|y|^2}{4t}}}
    \left(
      \log\gamma 
      +
      \log\frac{t}{|y|^2}
      +
      \log\frac{1}{t}
    \right)
    dy
    \\
    &
    =
    \log\gamma
    +
    \log\frac{1}{t}
    +c_2
    = 
    \log \frac{\gamma e^{c_2}}{2} + G^{-1}\left(\frac{t}{2}\right)
  \end{aligned}
\]
for $t>0$, where $c_2\in \mathbb{R}$ is a constant. 
Therefore, in the both cases, we see from \eqref{eq:4.3} that there exists $\gamma_*>0$ such that, 
for $\gamma > \gamma_*$, it holds 
\[
  \|e^{t\Delta}\varphi^{-1}(u_0)\|_{L^{\infty}}
  \ge G^{-1}\left(
    \frac{t}{2}
  \right)
\]
with small enough $t>0$.
Since $G$ is decreasing, by \eqref{eq:4.2} we obtain $F(\|e^{t\Delta}u_0\|_{L^{\infty}})\le \frac{t}{2}$ for $t>t_*$,
which contradicts Lemma~\ref{Lemma:2.2}.
Thus we complete the proof of Theorem~\ref{Theorem:1.3}. \qed

\section{Applications}
\label{section:5}


In this section, we consider the cases that
\begin{equation}
  \label{eq:5.1}
  f(u) = \exp\left( -\frac{1}{u} \right),
  \qquad u\in (0,1), 
\end{equation}
and 
\begin{equation}
  \label{eq:5.2}
  f(u) = u^p \left[ \log \left( e + \frac{1}{u} \right) \right]^{-r},
  \qquad u\in (0,1),
\end{equation}
where $1<p<1+2/N$ and $r>0$. 
The functions $f$ defined by \eqref{eq:5.1} and \eqref{eq:5.2} can be extended to functions whose domains are $(0,\infty)$
satisfying \eqref{eq:1.2} and \eqref{eq:1.3}.

We first focus on a nonlinear term satisfying \eqref{eq:1.2}, \eqref{eq:1.3} and \eqref{eq:5.1}.
\begin{Lemma}
  \label{Lemma:5.1}
  Let $N\ge 1$ and $f\in C^1([0,\infty)) \cap C^2((0,\infty))$ satisfy \eqref{eq:1.2}, \eqref{eq:1.3} and \eqref{eq:5.1}.
  Then
  \begin{equation*}
    f'(u)F(u) > 1
    = \lim_{u\to +0} f'(u)F(u)
  \end{equation*}
  for all sufficiently small $u>0$.
\end{Lemma}

\noindent
\textbf{Proof}.
By the l'H\^{o}pital rule
we have
\begin{equation*}
  \lim_{u\to +0} f'(u)F(u)
  =
  \lim_{u\to +0} \frac{(F(u))'}{(1/f'(u))'}
  =
  \lim_{u\to +0} \frac{f(u)f''(u)}{(f'(u))^2}
  = 1.
\end{equation*}
For $u\in (0,1)$, we have
\begin{align*}
  F(u)
  &=
  \int_u^1 (-s^2) \left( e^\frac{1}{s} \right)' \, ds + F(1)
  =
  u^2 e^\frac{1}{u} - e + 2 \int_u^1 s e^\frac{1}{s} \, ds + F(1)
  \\
  &\ge
  u^2 e^\frac{1}{u} - e + 2 \int_u^1 s \left( \frac{1}{2}\cdot \frac{1}{s^2} \right) \, ds + F(1)
  =
  u^2 e^\frac{1}{u} - e + \int_u^1 \frac{1}{s}\, ds + F(1)
  \\
  &=
  u^2 e^\frac{1}{u} - e + \log \frac{1}{u} + F(1)
  > u^2 e^\frac{1}{u}
\end{align*}
for all sufficiently small $u>0$.
Then we obtain
\begin{equation*}
  f'(u)F(u)
  >
  \frac{1}{u^2} e^{-\frac{1}{u}} \cdot u^2 e^\frac{1}{u}
  = 1
\end{equation*}
for all sufficiently small $u>0$.
Thus we complete the proof of Lemma~\ref{Lemma:5.1}.
\qed

\medskip

\begin{Theorem}
  \label{Theorem:5.1}
  Let $f\in C^1([0,\infty)) \cap C^2((0,\infty))$ satisfy \eqref{eq:1.2}, \eqref{eq:1.3} and \eqref{eq:5.1}. 
  For $c>0$, define 
  a nonnegative initial function 
  \begin{equation}
    \label{eq:5.3}
    u_{0}(x) = \frac{1}{\log \left[
      c^{-1} (|x|^2+ 1) \left\{
        \log (|x|+e)
      \right\}^2+1
    \right]}
  \end{equation}  in $\mathbb{R}^N$. 
  Then there exist positive constants $c_0\le c_1$ such that, 
  if $0<c<c_0$ then
  the solution $u$ of problem~\eqref{eq:1.1} exists globally in time; 
  if $c>c_1$ 
  then there cannot exist global in time solutions for problem~\eqref{eq:1.1}. 
\end{Theorem}

\noindent
\textbf{Proof of Theorem~\ref{Theorem:5.1}.}
We first consider the existence of a global solution. 
By Lemma~\ref{Lemma:5.1}, there exists $s_0>0$ such that 
$F(s) \ge 1/f'(s)$ for $0<s<s_0$. 
Let $c'<1$ be a positive constant satisfying
$\frac{1}{\log(1/c'+1)}<\min\{s_0,1\}$.
Take $c$ so that $0<c<c'$. 
Then there holds
$\|u_0\|_{L^{\infty}(\mathbb{R}^N)}<\min\{s_0,1\}$. 
Hence, it follows from
$F(s)\ge 1/f'(s) = s^2 e^{1/s}$ for $0<s<\min\{s_0,1\}$
that
\begin{equation}
  \label{eq:5.4}
  F(u_{0}(x))
  \ge \frac{1}{f'(u_0)}
=
\frac{1}{c}(|x|^2+1)
    \frac{\left\{\log (|x| + e)\right\}^2+1}
    {
    \left(
    \log 
    \left[
      c^{-1} (|x|^2+1) \left\{
        \log (|x|+e)
      \right\}^2+1
    \right]  
    \right)^2 
}
\end{equation}
for all $x\in \mathbb{R}^N$. 
We see from $\log (t+e) \le t+e$ for $t\ge 0$ and $c<1$ that 
\[
  \log \left[
    c^{-1} (|x|^2+1) \left\{
      \log (|x|+e)
    \right\}^2+1
  \right]  
  \le 
  \log \left[
    2c^{-1} (|x|+e)^4 
  \right] 
  = 
  \log \frac{2}{c} + 4 \log (|x|+e)
\]
for all $x\in\mathbb{R}^N$. 
Then, since $\frac{t}{\log\frac{2}{c}+4t}\ge \frac{1}{\log\frac{2}{c}+4}$ for $t\ge 1$,
by \eqref{eq:5.4} we obtain 
\[
  F(u_{0}(x))
  \ge
  \frac{1}{c}(|x|^2+1)
  \left(
    \frac{\log (|x| + e)}
    {
      \log \frac{2}{c} + 4 \log (|x|+e)
    }
  \right)^2 
\ge
\frac{1}{c\left(4+\log\frac{2}{c}\right)^2}(|x|^2+1)
\]
for all $x\in \mathbb{R}^N$.
Let $\gamma_*>0$ be the constant given in Theorem~\ref{Theorem:1.1}, 
and fix $\gamma\in (0,\gamma_*)$. 
Then there exists $c_0\in (0,c')$ such that
\[
  \frac{1}{c\left(4+\log\frac{2}{c}\right)^2}\ge \frac{1}{\gamma}
\]
for all $0<c<c_0$. 
This shows us that 
$u_0(x)\le F^{-1}\bigl[\gamma^{-1}(|x|^2+1)\bigr]$
for all $0<c<c_0$ and all $x\in \mathbb{R}^N$.
Then the solution constructed by Theorem~\ref{Theorem:1.1} is a global in time supersolution of problem \eqref{eq:1.1}, 
and hence there exists a global in time solution of \eqref{eq:1.1} with $u_0$ given in \eqref{eq:5.3}. 

We next consider the nonexistence of global solutions by applying
Theorem~\ref{Theorem:1.3}.
Let $\Gamma_*$ be the constant given in Theorem~\ref{Theorem:1.3}, 
and fix $c_1$ so that $c_1>2\Gamma_*$. 
By 
Lemma~\ref{Lemma:5.1} we can take $s_1>0$ such that 
$F(s)\le 2/f'(s)=2s^2e^{1/s}$ for all $0<s<s_1$. 
For any $c>c_1$ there exists $R>0$ such that 
$u_0(x) < s_1$
for all $x\in \mathbb{R}^N$ with $|x|\ge R$. 
Then we have 
$F(u_0)\le 2/f'(u_0)$ for all $|x|\ge R$. 
Therefore we obtain 
\[
  F(u_0)\le 
  \frac{2}{c}(|x|^2+1)
  \frac{\left\{\log (|x| + e)\right\}^2+1}
  {
    \left(
      \log 
      \left[
        c^{-1} (|x|^2+1) \left\{
          \log (|x|+e)
        \right\}^2+1
      \right]  
    \right)^2
  }
\]
for all $|x|\ge R$. 
This yields 
\[
  \liminf_{|x|\to \infty} \frac{|x|^2}{F(u_0(x))} 
  \ge 
  \frac{c}{2}
  > \frac{c_1}{2} > \Gamma_*
\]
for all $c>c_1$. 
Hence we can apply Theorem~\ref{Theorem:1.3}  
to prove the nonexistence of global solutions for problem~\eqref{eq:1.1}. 
This completes the proof of Theorem~\ref{Theorem:5.1}. 
\qed 

\medskip 

We next focus on 
a nonlinear term satisfying \eqref{eq:1.2}, \eqref{eq:1.3} and \eqref{eq:5.2}.

\begin{Theorem}
  \label{Theorem:5.2}
  For $1<p<1+2/N$ and $r>0$, 
  let $f\in C^1([0,\infty)) \cap C^2((0,\infty))$ satisfy \eqref{eq:1.2}, \eqref{eq:1.3} and \eqref{eq:5.2}. 
  For $c>0$, put 
  \begin{equation}
    \label{eq:5.5}
    u_{0}(x) = 
    c (|x|+1)^{-\frac{2}{p-1}}\left[
      \log \left(
        |x|+e
      \right)
    \right]^\frac{r}{p-1}
  \end{equation}
  in $\mathbb{R}^N$. 
  Then there exists positive constants $c_0\le c_1$ such that, 
  if $0<c<c_0$ then
  the solution $u$ of problem~\eqref{eq:1.1} exists globally in time; 
  if $c > c_1$ 
  then there cannot exist global in time solutions for problem~\eqref{eq:1.1}. 
\end{Theorem}

\noindent 
\textbf{Proof.}
Since $1<p<1+2/N$ and it follows from \eqref{eq:5.2} that 
\begin{gather}
  \label{eq:5.6}
  f'(s) = ps^{p-1}\left[\log \left(e+\frac{1}{s}\right)\right]^{-r} 
  + O\left(s^{p-1}\left[\log \left(e+\frac{1}{s}\right)\right]^{-r-1}\right), 
  \\ 
  \notag 
  f''(s) = p(p-1)s^{p-2}\left[\log \left(e+\frac{1}{s}\right)\right]^{-r} 
  + O\left(s^{p-2}\left[\log \left(e+\frac{1}{s}\right)\right]^{-r-1}\right), 
\end{gather}
as $s\to +0$,
we apply l'H\^{o}pital's rule to obtain 
\begin{equation}
  \label{eq:5.7}
  \lim_{s\to +0}f'(s)F(s) 
  =
  \lim_{s\to +0}\frac{\{f'(s)\}^2}{f(s)f''(s)}
  =
  \frac{p}{p-1}
  < 1 + \frac{N}{2}. 
\end{equation}
Hence $f$ satisfies the assumption in Theorem~\ref{Theorem:1.1}.
Let $\gamma_*>0$ be the constant given in Theorem~\ref{Theorem:1.1}, 
and fix $\gamma\in (0,\gamma_*)$. 
We show 
\begin{equation}
  \label{eq:5.8}
  u_0(x)\le F^{-1}\left[\gamma^{-1}(|x|^2+1)\right]
\end{equation}
for sufficiently small $c>0$. 
Then the solution constructed by Theorem~\ref{Theorem:1.1} is a global in time supersolution of problem \eqref{eq:1.1}, 
and hence there exists a global in time solution. 
By \eqref{eq:5.6} and \eqref{eq:5.7} we have 
$F(s) \ge 
\frac{3}{4}\frac{p}{p-1}\frac{1}{f'(s)}
\ge 
\frac{1}{2(p-1)}s^{-(p-1)}\left[\log \left(e+\frac{1}{s}\right)\right]^{r}$ for all $0<s<s_0$, 
where $s_0>0$ is a constant. 
Taking a sufficiently small $c'>0$, 
we may assume $\|u_0\|_{L^{\infty}(\mathbb{R}^N)}< s_0$ for all $0<c<c'$. 
Then we have 
\[
  F(u_0(x)) 
  \ge 
  \frac{1}{2(p-1)} c^{-(p-1)} (|x|+1)^{2}
  \left[
    \log \left(
      |x|+e
    \right)
  \right]^{-r}
  \left[ \log\left(e+c^{-1}\dfrac{(|x|+1)^\frac{2}{p-1}}{[\log(|x|+e)]^{\frac{r}{p-1}}}\right) \right]^r
\]
for all $x\in\mathbb{R}^N$. 
Since there exists a constant $m_1>0$ such that 
\[
  \frac{(|x|+1)^\frac{2}{p-1}}{[\log(|x|+e)]^{\frac{r}{p-1}}}
  \ge m_1 (|x|+e)^\frac{1}{p-1}
  \quad\mbox{in}\,\,\, \mathbb{R}^N, 
\]
for $c\in (0,\min\{c',m_1\})$ we obtain 
\begin{equation}
  \label{eq:5.9}
  \begin{aligned}
    F(u_0(x)) 
    & 
    \ge 
    \frac{1}{2(p-1)} c^{-(p-1)} (|x|+1)^{2}
    \left[
      \log \left(
        |x|+e
      \right)
    \right]^{-r}
    \left[ \log\left(e+c^{-1}m_1 (|x|+e)^\frac{1}{p-1} \right) \right]^r
    \\
    & 
    \ge 
    \frac{1}{2(p-1)} \left( \frac{1}{p-1} \right)^r c^{-(p-1)} (|x|^2+1) 
  \end{aligned}
\end{equation}
in $\mathbb{R}^N$. 
We now take a constant $c_0\in (0,\min\{c',m_1\})$ such that 
\begin{equation*}
  \frac{1}{2(p-1)} \left( \frac{1}{p-1} \right)^r c^{-(p-1)} 
  \ge 
  \frac{1}{\gamma}
\end{equation*}
for all $c\in (0,c_0)$. 
Then, by \eqref{eq:5.9} we obtain \eqref{eq:5.8} for all $c\in (0,c_0)$. 
Thus we obtain the existence of a global in time solution for problem~\eqref{eq:1.1} with $u_0$ given in \eqref{eq:5.5}. 

It remains to prove the nonexistence result for sufficiently large $c$.
Let $\Gamma_*$ be the constant given in Theorem~\ref{Theorem:1.3}, 
and fix $c_1$ so that
\begin{equation}
  \label{eq:5.10}
  \left( \frac{p-1}{2} \right)^{r+1} c_1^{p-1} > \Gamma_*. 
\end{equation} 
By \eqref{eq:5.6} and \eqref{eq:5.7} we have 
$F(s) \le 
\frac{3}{2}\frac{p}{p-1}\frac{1}{f'(s)}\le
\frac{2}{p-1}s^{-(p-1)}\left[\log \left(e+\frac{1}{s}\right)\right]^{r}$ for all $0<s<s_1$, 
where $s_1>0$ is a constant. 
For any $c>c_1$ there exists $R>0$ such that 
$0<u_0(x)<s_1$ for all $|x|\ge R$. 
Then it holds that
\[
  F(u_0(x))
  \le 
  \frac{2}{p-1} c^{-(p-1)} (|x|+1)^{2}
  \left[
    \log \left(
      |x|+e
    \right)
  \right]^{-r}
  \left[ \log\left(e+c^{-1}\dfrac{(|x|+1)^\frac{2}{p-1}}{[\log(|x|+e)]^{\frac{r}{p-1}}}\right) \right]^r
\]
for all $|x|\ge R$. 
This together with \eqref{eq:5.10} yields 
\[
  \liminf_{|x|\to \infty}\frac{|x|^2}{F(u_0(x))}
  \ge 
  \frac{p-1}{2} c^{p-1} \left( \frac{2}{p-1} \right)^{-r}
  =
  \left( \frac{p-1}{2} \right)^{r+1} c^{p-1} 
  > 
  \left( \frac{p-1}{2} \right)^{r+1} c_1^{p-1}
  > 
  \Gamma_*
\]
for all $c>c_1$. 
Thus we obtain the nonexistence of global in time solutions for problem~\eqref{eq:1.1}. 
This completes the proof of Theorem~\ref{Theorem:5.2}.
\qed 

\medskip 

\noindent
{\bf Acknowledgements.}
The first author was supported partially by the Grant-in-Aid for Early-Career Scientists (No.~19K14569).
The second author
was supported in part
by the Grant-in-Aid for Challenging Exploratory Research(No.~19KK0349)
and the Grant-in-Aid for Fostering Joint International Research(A)
(No.~21K18582)
from Japan Society for the Promotion of Science.



\end{document}